\documentclass[11pt]{article}
\usepackage{amsfonts}
\usepackage{bbm}
\usepackage{mathrsfs}
\usepackage{amssymb}
\usepackage{url}

\title{The Standard Quantized Matrix Algebra $M_q(n)$\\ is A Solvable Polynomial Algebra}

\author{Rabigul Tuniyaz\thanks{Project supported by
the National Natural Science Foundation of China (1186106).}\\
{\small Department of Mathematics, School of  Science}\\
{\small Xinjiang Institute of Science and Technology}\\
{\small  Akesu, 843100, Xinjiang, China}\\}

\date{}

\begin{document}
\maketitle
\begin{center}
\begin{minipage}{120mm}
{\small {\bf Abstract.} Let $M_q(n)$ be the standard quantized matrix algebra, introduced by Faddeev, Reshetikhin, and Takhtajan. It is  shown, by constructing an appropriate monomial ordering $\prec$ on its PBW $K$-basis ${\cal B}$ , that $M_q(n)$ is a solvable polynomial algebra. Consequently, further structural properties of $M_q(n)$  and their modules may be established and realized in a constructive-computational way.}
\end{minipage}\end{center} {\parindent=0pt\par

}
\vskip -.5truecm

\renewcommand{\thefootnote}{\fnsymbol{footnote}}
\let\footnote\relax\footnotetext{E-mail: rabigul802@sina.com}
\let\footnote\relax

\def\NZ{\mathbb{N}}
\def\QED{\hfill{$\Box$}}
\def \r{\rightarrow}
\def\mapright#1#2{\smash{\mathop{\longrightarrow}\limits^{#1}_{#2}}}

\def\v5{\vskip .5truecm}
\def\OV#1{\overline {#1}}
\def\hang{\hangindent\parindent}
\def\textindent#1{\indent\llap{#1\enspace}\ignorespaces}
\def\item{\par\hang\textindent}

\def\LH{{\bf LH}}\def\LM{{\bf LM}}\def\LT{{\bf
LT}}\def\KX{K\langle X\rangle} \def\KS{K\langle X\rangle}
\def\B{{\cal B}} \def\LC{{\bf LC}} \def\G{{\cal G}} \def\FRAC#1#2{\displaystyle{\frac{#1}{#2}}}
\def\SUM^#1_#2{\displaystyle{\sum^{#1}_{#2}}}
\def\T{\widetilde}\def\KD{K\langle D\rangle}
\def\PRC{\prec_{d\textrm{\tiny -}lex}}
\def\KZ{K\langle Z\rangle}

\section*{1. Introduction}
Let $K$ be a field of characteristic 0. The standard quantized matrix algebra $M_q(n)$,   introduced in [FRT], has been widely studied and generalized in different contexts, for instance, [JP], [JZ1, 2],  and [JJJ]. In [Tu], it was  shown  explicitly that the defining relations of $M_q(n)$ form a Gr\"obner-Shirshov basis $\G$, and by means of $\G$, several structural properties of $M_q(n)$ were derived, for instance, $M_q(n)$ has a PBW $K$-basis, is  of Hilbert series $\frac{1}{(1-t)^{n^2}}$, of Gelfand-Kirillov dimension $n^2$, of global homological dimension $n^2$,  is a classical Koszul algebra,  and has the elimination property for one-sided ideals $L$ with GK.dim$M_q(n)/L<$ GK.dim$M_q(n)$. The main purpose of this note is to show that $M_q(n)$ is a solvable polynomial algebra in the sense of [K-RW] (Section 2, Theorem 2.3), which means that $M_q(n)$ has an algorithmic Gr\"obner basis theory for both two-sided and one-sided ideals. To demonstrate how  Theorem 2.3 may bring some perspective of establishing  and realizing further structural properties of $M_q(n)$ and their modules in a constructive-computational way, in Section 3 we specify several applications of Theorem 2.3. \v5

Throughout this note, $K$ denotes a field of characteristic 0, $K^*=K-\{0\}$, and all $K$-algebras considered are associative with multiplicative identity 1. If $S$ is a nonempty subset of an algebra $A$, then we write $\langle S\rangle$ for the two-sided ideal of $A$ generated by $S$.

\section*{2. $M_q(n)$ is a solvable polynomial algebra}
We start  by recalling from ([K-RW], [Li1, 6]) the following definitions and notations. Suppose that a finitely generated  $K$-algebra $A=K[a_1,\ldots ,a_n]$ has the PBW $K$-basis $\B =\{ a^{\alpha}=a_{1}^{\alpha_1}\cdots
a_{n}^{\alpha_n}~|~\alpha =(\alpha_1,\ldots ,\alpha_n)\in\NZ^n\}$, and that $\prec$ is a total ordering on $\B$. Then every nonzero element $f\in A$ has a unique expression
$$\begin{array}{rcl} f&=&\lambda_1a^{\alpha (1)}+\lambda_2a^{\alpha (2)}+\cdots +\lambda_ma^{\alpha (m)},\\
&{~}&\hbox{such that}~a^{\alpha (1)}\prec a^{\alpha
(2)}\prec\cdots \prec a^{\alpha (m)},\\
&{~}&\hbox{where}~ \lambda_j\in K^*,~a^{\alpha
(j)}=a_1^{\alpha_{1j}}a_2^{\alpha_{2j}}\cdots a_n^{\alpha_{nj}}\in\B
,~1\le j\le m.
\end{array}$$
Since elements of $\B$ are conventionally called {\it monomials}\index{monomial}, the {\it leading monomial of $f$} is defined as $\LM
(f)=a^{\alpha (m)}$, the {\it leading coefficient of $f$} is defined
as $\LC (f)=\lambda_m$, and the {\it leading term of $f$} is defined
as $\LT (f)=\lambda_ma^{\alpha (m)}$.\v5

{\bf Definition 2.1}  Suppose that the $K$-algebra
$A=K[a_1,\ldots ,a_n]$ has the PBW basis $\B$. If $\prec$ is a
total ordering on $\B$ that satisfies the following three
conditions:{\parindent=1.35truecm\par

\item{(1)} $\prec$ is a well-ordering (i.e., every nonempty subset of $\B$ has a minimal element);\par

\item{(2)} For $a^{\gamma},a^{\alpha},a^{\beta},a^{\eta}\in\B$, if $a^{\gamma}\ne 1$, $a^{\beta}\ne
a^{\gamma}$, and $a^{\gamma}=\LM (a^{\alpha}a^{\beta}a^{\eta})$,
then $a^{\beta}\prec a^{\gamma}$ (thereby $1\prec a^{\gamma}$ for
all $a^{\gamma}\ne 1$);\par

\item{(3)} For $a^{\gamma},a^{\alpha},a^{\beta}, a^{\eta}\in\B$, if
$a^{\alpha}\prec a^{\beta}$, $\LM (a^{\gamma}a^{\alpha}a^{\eta})\ne
0$, and $\LM (a^{\gamma}a^{\beta}a^{\eta})\not\in \{ 0,1\}$, then
$\LM (a^{\gamma}a^{\alpha}a^{\eta})\prec\LM
(a^{\gamma}a^{\beta}a^{\eta})$,\par}{\parindent=0pt
then $\prec$ is called a {\it monomial ordering} on $\B$ (or a
monomial ordering on $A$).} \v5

{\bf Definition 2.2} A finitely generated $K$-algebra $A=K[a_1,\ldots ,a_n]$
is called a {\it solvable polynomial algebra} if $A$ has the PBW $K$-basis $\B =\{
a^{\alpha}=a_1^{\alpha_1}\cdots a_n^{\alpha_n}~|~\alpha
=(\alpha_1,\ldots ,\alpha_n)\in\NZ^n\}$ and a monomial ordering $\prec$ on $\B$, such that
for $\lambda_{ji}\in K^*$ and  $f_{ji}\in A$,
$$\begin{array}{l} a_ja_i=\lambda_{ji}a_ia_j+f_{ji},~1\le i<j\le n,\\
\LM (f_{ji})\prec a_ia_j~\hbox{whenever}~f_{ji}\ne 0.\end{array}$$\par

Now, we aim to prove the following result.\v5

{\bf Theorem 2.3}  Let $M_q(n)$ be the standard quantized matrix algebra over a field $K$, in the sense of [FRT]. Then $M_q(n)$ is a solvable polynomial algebra in the sense of Definition 2.2. \vskip 6pt

{\bf Proof} Let $I(n)=\{(i,j)~|~i,j=1,2,\cdots,n\}$ with $n\ge 2$. Recall from [FRT] that $M_q(n)$ is the associative $K$-algebra generated by the set of $n^2$ generators
$Z=\{ z_{ij}~|~(i,j)\in I(n)\}$  subject to the relations:
$$\begin{array}{ll}
\hbox{R}_1:\quad z_{ij}z_{ik}=qz_{ik}z_{ij},&\hbox{if}~j<k ,\\
\hbox{R}_2:\quad z_{ij}z_{kj}=qz_{kj}z_{ij},&\hbox{if}~i<k ,\\
\hbox{R}_3:\quad z_{ij}z_{st}=z_{st}z_{ij},&\hbox{if}~i<s,\;t<j,\\
\hbox{R}_4:\quad z_{ij}z_{st}=z_{st}z_{ij}+(q-q^{-1})z_{it}z_{sj},&\hbox{if}~i<s,~j<t,
 \end{array}$$
where $i,j,k,s,t=1,2,...,n$ and $q\in K^*$ is the quantum parameter. By [Rab, Corollary 3.1], $M_q(n)$ has the  PBW $K$-basis
$$\B =\left\{\left. z^{k_{nn}}_{nn}z^{k_{nn-1}}_{nn-1}\cdots z^{k_{n1}}_{n_1}z^{k_{n-1n}}_{n-1n}\cdots z^{k_{n-11}}_{n-11}\cdots z^{k_{1n}}_{1n}\cdots z^{k_{11}}_{11}~\right |~k_{ij}\in \NZ,(i, j)\in I(n)\right\}.$$ \par

We now start on constructing a monomial ordering on $\B$ such that all conditions of Definition 2.1 and Definition 2.2 are satisfied. In doing so, we first rewrite $\B$ as
$$\B =\left\{ 1,~z_{i_1j_1}z_{i_2j_2}\cdots z_{i_kj_k}~\left |~
\begin{array}{l} (i_q,j_q)\in I(n),~k\ge 1, \\
(i_1,j_1)\ge (i_2,j_2)\ge\cdots\ge (i_k,j_k)\end{array}\right.\right\} ,$$
where
$$(i_{\ell},j_{\ell})<(i_t,j_t)\Leftrightarrow\left\{
\begin{array}{l} i_{\ell}<i_t,\\
\hbox{or}~i_{\ell}=i_t~\hbox{and}~j_{\ell}<j_t.\end{array}\right.$$
Then, we define the ordering $\prec$ on the set $Z$ of generators: for $z_{kl}$, $z_{ij}\in Z$,
$$z_{kl}\prec z_{ij}\Leftrightarrow\left\{\begin{array}{l}
k<i,\\
\hbox{or}~k=i~\hbox{and}~l <j.\end{array}\right.$$
and extend this ordering to $\B$:
$$1\prec u~\hbox{for all}~u=z_{k_1l_1}z_{k_2l_2}\cdots z_{k_rl_r}\in\B -\{ 1\},$$
and for $u=z_{k_1l_1}z_{k_2l_2}\cdots z_{k_rl_r},$, $v=z_{i_1j_1}z_{i_2j_2}\cdots z_{i_qj_q}\in\B -\{ 1\}$,
$$u\prec v\Leftrightarrow\left\{\begin{array}{l} r<q~\hbox{and}~ z_{k_1l_1}=z_{i_1j_1},z_{k_2l_2}=z_{i_2j_2},\ldots ,z_{k_rl_r}=z_{i_rj_r},\\
\hbox{or there exists am}~m,~1\le m\le r,~\hbox{such that}\\
z_{k_1l_1}=z_{i_1j_1},z_{k_2l_2}=z_{i_2j_2},\ldots ,z_{k_{m-1}l_{m-1}}=z_{i_{m-1}j_{m-1}},~
\hbox{but}~z_{k_ml_m}\prec z_{i_mj_m}.\end{array}\right.$$
It is straightforward to check that $\prec$ is reflexive, antisymmetrical, transitive, and any two elements $u,v\in\B$ are comparable, thereby $\prec$ is a
total ordering on $\B$. Also since $I(n)$ is a finite set, it can be directly verified  that $\prec$ satisfies the descending chain condition on $\B$, namely $\prec$ is a well-ordering on $\B$. \par
It remains to show that $\prec$ satisfies the conditions (2) and (3) of Definition 2.1, and that with respect to $\prec$ on $\B$, the relations $\hbox{R}_1$, $\hbox{R}_2$, $\hbox{R}_3$, $\hbox{R}_4$ satisfied by generators of $M_q(n)$ have the property required by Definition 2.2. To this end, we first observe that in the relations $\hbox{R}_1$, $\hbox{R}_2$, $\hbox{R}_3$, and $\hbox{R}_4$, the monomials $z_{ik}z_{ij}, z_{kj}z_{ij}, z_{st}z_{ij}\in \B$. Next, let $z_{ij}, z_{st}, z_{pq}\in Z$, and suppose that $z_{st}\prec z_{pq}$. If $(i,j)$ and $(s,t)$ are such that $i<s$ and $j<t$, then the relation $\hbox{R}_3$ gives rise to
$$z_{it}z_{sj}=z_{sj}z_{it}\in\B ,$$
thereby the relation $\hbox{R}_4$ is turned into
$$\begin{array}{rcl} z_{ij}z_{st}&=&z_{st}z_{ij}+(q-q^{-1})z_{sj}z_{it}\\
&{~}&\hbox{with}~z_{st}z_{ij},z_{sj}z_{it}\in\B\\
&{~}&\hbox{and}~z_{sj}z_{it}\prec z_{st}z_{ij}=\LM (z_{ij}z_{st}).\end{array}$$
On the other hand, if $(i,j)$ and $(p,q)$ are such that $i<p$ and $j<q$, then the the relation $\hbox{R}_3$ gives rise to
$$z_{iq}z_{pj}=z_{pj}z_{iq}\in\B ,$$
thereby the relation $\hbox{R}_4$ is turned into
$$\begin{array}{rcl} z_{ij}z_{pq}&=&z_{pq}z_{ij}+(q-q^{-1})z_{pj}z_{iq}\\
&{~}&\hbox{with}~z_{pq}z_{ij},z_{pj}z_{iq}\in\B\\
&{~}&\hbox{and}~z_{pj}z_{iq}\prec z_{pq}z_{ij}=\LM (z_{ij}z_{pq}).\end{array}$$
Thus, we have shown that if
$$\begin{array}{l} (i,j),(s,t)\in I(n)~\hbox{such that}~i<s,j<t,\\
(i,j),(p,q)\in I(n)~\hbox{such that}~i<p,j<q,\end{array}$$
then
$$\begin{array}{l} z_{st}\prec z_{pq}~\hbox{implies}~\LM (z_{ij}z_{st})=z_{st}z_{ij}\prec z_{pq}z_{ij}=\LM (z_{ij}z_{pq}),\\
\hbox{and the generating relations of}~M_q(n)~\hbox{determined by}~\hbox{R}_4\\
\hbox{have the property required by Definition 2.2}.\end{array}\eqno{(1)}$$
Similarly in the case that
$$\begin{array}{l} (s,t), (i,j)\in I(n)~\hbox{such that}~s<i, t<j,\\
(p, q), (i,j)\in I(n)~\hbox{such that}~p<i, q<j,\end{array}$$
with the aid of $\hbox{R}_3$ we have
$$\begin{array}{l} z_{st}\prec z_{pq}~\hbox{implies}~\LM (z_{st}z_{ij})=z_{ij}z_{st}\prec z_{ij}z_{pq}=\LM (z_{pq}z_{ij}),\\
\hbox{and the generating relations of}~M_q(n)~\hbox{determined by}~\hbox{R}_4\\
\hbox{have the property required by Definition 2.2}.\end{array}\eqno{(2)}$$
At this stage, bearing in mind the relations $\hbox{R}_1$, $\hbox{R}_2$, $\hbox{R}_3$, $\hbox{R}_4$, and the assertions (1) and (2) derived above, we may conclude that
$$\begin{array}{l} \hbox{for any}~
z_{ij}, z_{st}, z_{pq}\in Z,~\hbox{if} z_{st}\prec z_{pq},~\hbox{then}\\
\LM (z_{ij}z_{st})\prec \LM (z_{ij}z_{pq}), \LM (z_{st}z_{ij})\prec \LM (z_{pq}z_{ij}),
~\hbox{and}\\
\hbox{the generating relations of}~M_q(n)~\hbox{determined by}~\hbox{R}_1, \hbox{R}_2, \hbox{R}_3,  \hbox{and R}_4\\
\hbox{have the property required by Definition 2.2}.\end{array}\eqno{(3)}$$\par

Finally, by means of (1), (2), and (3) presented above, it is straightforward to check that the conditions (2) and (3) of Definition 2.1 are satisfied by $\prec$, thereby $\prec$ is a monomial ordering on $\B$, and consequently $M_q(n)$ is a solvable polynomial algebra in the sense of Definition 2.2, as desired.\QED\v5

\section*{3. Several applications of Theorem 2.3}

From Theorem 2.3 obtained in the last section we have known that the standard quantized matrix algebra $M_q(n)$ is a solvable polynomial algebra in the sense of [K-RW]. Thus,
it is well known that every (two-sided, respectively one-sided) ideal of a solvable polynomial algebra $A$ and every submodule of a free (left) $A$-module has a finite Gr\"obner basis with respect to a given monomial ordering, in particular, for one-sided ideals and submodules of free (left) modules there is a noncommutative Buchberger Algorithm which, nowadays, has been successfully implemented in the computer algebra system \textsf{Plural} [LS].  At this point, we specify several applications of Theorem 2.3 in this section, so as to show how  Theorem 2.3 may bring some perspective of establishing and realizing  further structural properties of $M_q(n)$ and their modules in a constructive-computational way. For more details on the basic constructive-computational theory and methods for solvable polynomial algebras and their modules, one is referred to [Li6]. All notations used in Section 2 are maintained. Moreover, modules over $M_q(n)$ are meant {\it left $M_q(n)$-modules}.\v5

{\bf Theorem 3.1} Let $M_q(n)$ be the standard quantized matrix algebra. Then the following statements hold.\par

(i) $M_q(n)$ is a Noetherian domain.\par

(ii) Let $L$ be a nonzero left ideal of $M_q(n)$, and $M_q(n)/L$ the left $M_q(n)$-module. Considering Gelfand-Kirillov dimesion, we have GK.dim$M_q(n)/L<$ GK.dim$M_q(n)=n^2$, and there is an algorithm for computing GK.dim$M_q(n)/L$.\par

(iii)  Let $M$ be a finitely generated  $M_q(n)$-module. Then a finite free resolution of $M$ can be algorithmically constructed, and the projective dimension of $M$ can be algorithmically computed.\par

(iv) Let $M$ be a finitely generated graded $M_q(n)$-module (note that $M_q(n)$ is an $\NZ$-graded algebra in which each generator has the degree 1).  Then a minimal homogeneous generating set of $M$ can be algorithmically computed, and a minimal finite graded free resolution of $M$ can be algorithmically constructed.\vskip 6pt

{\bf Proof} (i) Though the property that $M_q(n)$ is a Noetherian domain may be (or may have been) established  in some other ways,  here we emphasize that this property may follow immediately from Theorem 2.3. More precisely, that $M_q(n)$ has no divisors of zero follows from the fact that $\LM (fg)=\LM (f)\LM (g)$ for all nonzero $f,g\in M_q(n)$, and that the Noetherianess of $M_q(n)$ follows from the fact that every nonzero one-sided ideal has a finite Gr\"obner basis (see [K-RW]). \par

(ii) That Gk.dim$D_q(n)=n^2$ follows from [Rab]. Since $M_q(n)$ is a (quadratic) solvable polynomial algebra by Theorem 2.3, it follows from [Li1, CH.V] that GK.dim$M_q(n)/L<n^2$ (this may also follow from classical Gelfand-Kirillov dimension theory [KL], for $M_q(n)$ is now a Noetherian domain), and that there is an algorithm for computing GK.dim$M_q(n)/L$.\par

(iii) This follows from [Li6, Ch.3].\par

(iv) This follows from [Li6, Ch.4].\QED\v5

In [Tu], it was shown that every left ideal $L$ of $M_q(n)$  with Gelfand-Kirillov dimension GK.dim$M_q(n)/L<$ GK.dim$M_q(n)$ has the elimination property in the sense of [Li5]. Now that we know that $M_q(n)$ is a solvable polynomial algebra, [Li6, Section 1.6] tells us that this elimination property can be  strengthened, and may be realized in a computational way.  To see this clearly, let us first recall the Elimination Lemma given in [8].  Let $A=K[a_1,\ldots ,a_n]$ be a  finitely generated $K$-algebra with
the PBW basis $\B =\{a^{\alpha}=a_1^{\alpha_1}\cdots
a_n^{\alpha_n}~|~\alpha =(\alpha_1,\ldots
,\alpha_n)\in\NZ^n\}$ and, for a subset $U=\{ a_{i_1},...,a_{i_r}\}\subset\{
a_1,...,a_n\}$ with $i_1<i_2<\cdots <i_r$,  let
$$S=\left\{ a_{i_1}^{\alpha_1}\cdots a_{i_r}^{\alpha_r}~\Big |~
(\alpha_1,...,\alpha_r)\in\NZ^r\right\},\quad
V(S)=K\hbox{-span}S.$$

{\bf Lemma 3.2} {[8, Lemma 3.1]} With notation as fixed
above, let $L$ be a nonzero left ideal of $A$ and $A/L$ the left $A$-module defined by $L$. If there is a subset  $U=\{ a_{i_1},\ldots ,a_{i_r}\}
\subset\{a_1,\ldots ,a_n\}$ with $i_1<i_2<\cdots <i_r$, such that $V(S)\cap L=\{ 0\}$, then $$\hbox{GK.dim}(A/L)\ge r.$$  Consequently, if  $A/L$ has finite GK dimension $\hbox{GK.dim}(A/L)=d<n$ ($=$
the number of generators of $A$), then  $$V(S)\cap L\ne
\{ 0\}$$ holds true for every subset $U=\{
a_{i_1},...,a_{i_{d+1}}\}\subset$ $\{ a_1,...,a_n\}$ with
$i_1<i_2<\cdots <i_{d+1}$, in particular, for every $U=\{
a_1,\ldots a_s\}$ with $d+1\le s\le n-1$, we have $V(S)\cap
L\ne \{ 0\}$.\par\QED\v5

For convenience of stating the next theorem, let us write the set of generators of $M_q(n)$ as $Z=\{ z_1,z_2,\ldots z_{n^2}\}$, i.e., $M_q(n)=K[z_1,z_2,\ldots ,z_{n^2}]$.\v5

{\bf Theorem 3.3} With notation as fixed above, Let $L$ be a nonzero left ideal of $M_q(n)$. Then the following two statements hold. \par

(i) GK.dim$M_q(n)/L< n^2=$ GK.dim$M_q(n)$. If GK.dim$M_q(n)/L=d$, then
$$V(T)\cap L\ne\{ 0\}$$ holds true for every subset $U=\{
z_{i_1},z_{i_2},...,z_{d+1}\}\subset Z$ with
$i_1<i_2<\cdots <i_{d+1}$, in particular, for every $U=\{
x_1,x_2\ldots x_s\}$ with $d+1\le s\le n^2 -1$, we have $V(T)\cap
L\ne \{ 0\}$.\par

(ii) Without exactly knowing the numerical value GK.dim$M_q(n)/L$, the elimination property for a left ideal  $L=\sum_{i=1}^mM_q(n)\xi_i$ of $M_q(n)$ can be realized in a computational way, as follows:\par

Let $\prec$ be the monomial ordering on the PBW basis $\B$ of $M_q(n)$ as constructed in the proof of Theorem 2.3, and let $V(T)$ be as in (i). Then, employing an elimination ordering  $\lessdot$ with respect to $V(T)$ (which can always  be constructed if the existing monomial ordering on $\B$ is not an elimination ordering, see [Li6, Proposition 1.6.3]),  a Gr\"obner basis $\G$ of $L$ can be produced by running the noncommutative Buchberger algorithm for solvable polynomial algebras, such that
$$L\cap V(T)\ne \{ 0\} \Leftrightarrow \G\cap V(T)\ne \emptyset .$$
\vskip 6pt

{\bf Proof} (i) Since $M_q(n)$ has the PBW basis $\B$, GK.dim$M_q(n)=n^2$ by [Rab], and GK.dim$M_q(n)/L<n^2$ by Theorem 3.1(ii), the desired elimination property follows from  Lemma 3.2 mentioned above.\par

(ii) This follows from [Li6, Corollary 1.6.5].\QED\v5

{\bf Remark} Since $M_q(n)$ is now a solvable polynomial algebra, if  $F=\oplus_{i=1}^sM_q(n)e_i$ is a free (left) $M_q(n)$-module of finite rank, then a similar (even much stronger) result of Theorem 3.3 holds true for any finitely generated submodule $N=\sum_{i=1}^mM_q(n)\xi_i$ of $F$. The reader is referred to [Li6, Section 2.4] for a detailed argumentation.\v5

\centerline{Refeerence}{\parindent=1.47truecm\par

\item{[Bok]} L. Bokut et al., {\it Gr\"obner--Shirshov Bases:
Normal Forms, Combinatorial and Decision Problems in Algebra}. World Scientific Publishing,
2020. \url{https://doi.org/10.1142/9287}\par

\item{[FRT]} L. D. Faddeev, N. Yu. Reshetikhin,  L. A. Takhtajan, Quantization of Lie groups and Lie algebras. {\it Algebraic Analysis}, Academic Press  (1988), 129--140.\par

\item{[JP]} H. P. Jakobsen, C. Pagani, Quantized matrix algebras and quantum seeds. {\it Linear and Multilinear Algebra}, 2014. DOI: 10.1080/03081087.2014.898297

\item{[JZ1]} H. P. Jakobsen, H. Zhang, The center of the quantized matrix algebra. {\it J. Algebra}, (196)(1997), 458--474.

\item{[JZ2]} H. P. Jakobsen and H. Zhang , A class of quadratic matrix algebras arising from the quantized enveloping algebra $U_q(A_{2n-1})$. {\it J. Math. Phys}., (41)(2000), 2310--2336.

\item{[JJJ]} H. P. Jakobsen, S. J\"ondrup, A. Jensen, Quadratic algebras of type AIII.III.  In: {\it Tsinghua Science $\&$ Technology}, (3)(1998), 1209--1212 .

\item{[KL]} G.R. Krause, T.H. Lenagan, {\it Growth of Algebras and
Gelfand-Kirillov Dimension}. Graduate Studies in Mathematics.
American Mathematical Society, 1991.

\item{[K-RW]} A. Kandri-Rody, V. Weispfenning, Non-commutative
Gr\"obner bases in algebras of solvable type. {\it J. Symbolic
Comput.}, 9(1990), 1--26. Also available as: Technical Report University of Passau,
MIP-8807, March 1988.

\item{[Li1]} H. Li, {\it Noncommutative Gr\"obner Bases and Filtered-graded Transfer}.
Lecture Notes in Mathematics, Vol. 1795, Springer, 2002.

\item{[Li2]} H. Li, $\Gamma$-leading ~homogeneous~ algebras~ and
Gr\"obner bases. In: {\it Recent Developments in Algebra and
Related Areas} (F. Li and C. Dong eds.), Advanced Lectures in
Mathematics, Vol. 8, International Press \& Higher Education Press,
Boston-Beijing, 2009, 155 -- 200. Also available at: arXiv:math/0609583 [math.RA].

\item{[Li3]} H. Li, {\it Gr\"obner Bases in Ring Theory}. World Scientific Publishing Co., 2011. \url{https://doi.org/10.1142/8223}\par

\item{[Li4]} H. Li, A note on solvable polynomial algebras. {\it Computer Science Journal of Moldova}, vol.22, no.1(64), 2014, 99--109.  arXiv:1212.5988 [math.RA]

\item{[Li5]} H. Li, An elimination lemma for algebras with PBW bases. {\it Communications in
Algebra}, 46(8)(2018), 3520--3532.\par

\item{[Li6]} H. Li, {\it Noncommutative polynomial algebras of solvable type and their modules: Basic constructive-computational theory and methods}. Chapman and Hall/CRC Press, 2021.

\item{[LS]} V. Levandovskyy, H. Sch\"onemann, Plural: a computer algebra system for noncommutative polynomial algebras. In: {\it Proc. Symbolic and Algebraic Computation}, International Symposium ISSAC 2003, Philadelphia, USA, 176--183, 2003.

\item{[Tu]} R. Tuniyaz,  Some Structural Properties of the Standard Quantized Matrix Algebra Mq(n).	arXiv:2112.13628[math.RA]

\end{document}